\input amstex
\magnification=1200
\documentstyle{amsppt}
\NoBlackBoxes
\widestnumber\key{[Most]2]} 
\baselineskip 20pt

\hsize 6.25truein
\vsize = 8.75truein


\catcode`\@=11
\def\logo@{}
\catcode`\@=13

\def\cplus{\hbox{$\subset${\raise0.3ex\hbox{\kern -0.55em
${\scriptscriptstyle +}$}}\ }}

\def\clplus{\hbox{$\subset${\raise0.3ex\hbox{\kern -0.55em ${\scriptscriptstyle +}$}}\ }}
\def\crplus{\hbox{$\supset${\raise1.05pt\hbox{\kern -0.55em ${\scriptscriptstyle +}$}}\ }}

\def\N{\bold N}

\def\p{\bold p}
\def\C{\bold C}
\def\Z{\bold Z}
\def\R{\bold R}

\def\id{\text{\rm id}}

\def\supp{\text{\rm supp}}
\def\Re{\text{\rm Re}}
\def\Im{\text{\rm Im}}
\def\Ext{\text{\rm Ext}}

\def\p{\bold p}
\def\C{\bold C}
\def\Z{\bold Z}
\def\R{\bold R}

\def\id{\text{\rm id}}

\def\ch{\text{\rm ch}}
\def\pro{\text{\rm pro}}
\def\Hom{\text{\rm Hom}}
\def\m{\frak m}
\def\n{\frak n}
\def\p{\frak p}
\def\g{\frak g}
\def\h{\frak h}
\def\b{\frak b}
\def\k{\frak k}
\def\t{\frak t}

\centerline{\qquad\qquad\qquad\qquad\qquad\it To the memory of Armand Borel}
\vskip .15in
\topmatter
\title GENERALIZED HARISH-CHANDRA MODULES WITH GENERIC MINIMAL $\frak k$-type \endtitle
\author Ivan Penkov and Gregg Zuckerman \endauthor
\affil  \endaffil	
\abstract  We make a first step towards a classification of simple generalized Harish-Chandra modules which are not Harish-Chandra modules or weight modules of finite type.  For an arbitrary algebraic reductive pair of complex Lie algebras $(\g,\k)$, we construct, via cohomological induction, the fundamental series $F^\cdot (\p,E)$ of generalized Harish-Chandra modules.  We then use $F^\cdot (\p,E)$ to characterize any simple generalized Harish-Chandra module with generic minimal $\k$-type.  More precisely, we prove that any such simple $(\g,\k)$-module of finite type arises as the unique simple submodule of an appropriate fundamental series module $F^s(\p,E)$ in the middle dimension $s$. Under the stronger assumption that $\k$ contains a semisimple regular element of $\g$, we prove that any simple $(\g,\k)$-module with generic minimal $\k$-type is necessarily of finite type, and hence obtain a reconstruction theorem for a class of simple $(\g,\k)$-modules which can a priori have infinite type.  We also obtain generic general versions of some classical theorems of Harish-Chandra, such as the Harish-Chandra admissibility theorem.  The paper is concluded by examples, in particular we compute the genericity condition on a $\k$-type for any pair $(\g,\k)$ with $\k\simeq s\ell (2)$.
\endabstract
\endtopmatter

\head Introduction \endhead

The goal of the present paper is to make a first step towards a classification of simple generalized Harish-Chandra modules which are not Harish-Chandra modules or weight modules of finite type.  This work is part of the program of study of generalized Harish-Chandra modules laid out in \cite{PZ}.  Let $\frak g$ be a semisimple Lie algebra.  A {\it simple generalized Harish-Chandra module} is by definition a simple $\frak g$-module with locally finite action of a reductive in $\frak g$ subalgebra $\k\subset \frak g$ and with finite $\k$-multiplicities.  In the classical case of Harish-Chandra modules, the pair $(\frak g,\frak k)$ is in addition assumed to be symmetric.  In a recent joint paper with V. Serganova \cite{PSZ}, we have constructed new families of generalized Harish-Chandra modules; however, no general classification is known beyond the case when the pair $(\frak g, \frak k)$ is symmetric and the case when $\frak k$ is a Cartan subalgebra of $\frak g$.  The first case is settled in the well-known work of R. Langlands \cite{L2}, A. Knapp and the second named author \cite{KZ}, D. Vogan and the second named author \cite{V2}, \cite{Z}, A. Beilinson - J. Bernstein \cite{BB} and I. Mirkovic  \cite{Mi}; the second case is settled in a more recent breakthrough by O. Mathieu \cite{M}.  Some low rank cases of certain special non-symmetric pairs  $(\frak g,\frak k)$ (where $\k$ is not a Cartan subalgebra) have been settled by G. Savin \cite{Sa}.

In this paper, we consider simple generalized Harish-Chandra modules which have a generic minimal $\frak k$-type for some arbitrary fixed reductive pair $(\frak g,\frak k)$ (the precise definitions see in Section 1 below).  One of  our main results is the construction of a series of $(\frak g,\frak k)$-modules, which we call the fundamental series (it generalizes the fundamental series of Harish-Chandra modules), and furthermore the theorem that any simple generalized Harish-Chandra module with generic minimal $\frak k$-type is a submodule of the fundamental series.  We refer to the latter result as the first reconstruction theorem for generalized Harish-Chandra modules.  This theorem is based on new results on the $\n$-cohomology of a simple generalized Harish-Chandra module and on the vanishing of cohomological induction except in the middle dimension (see Theorem 1 and 2 in Section 1 below).  Under a stronger assumption on the pair $(\frak g,\frak k)$, we also prove a reconstruction theorem for simple $(\frak g,\frak k)$-modules which may a priori have infinite type (we refer to it as the second reconstruction theorem) and in particular a new generic general version of Harish-Chandra's admissibility theorem, see Theorem 4 and Corollary 2 in Section 1.

Here is a very brief historical perspective.  In the 1950's, the classical Borel-Weil-Bott theorem opened a new era in representation theory, relating the simple finite dimensional representations of a semisimple Lie algebra $\frak g$ with the, possibly higher, cohomology of homogeneous line bundles on the flag variety of a complex algebraic group $G$ with Lie algebra $\frak g$, \cite{S}, \cite{Bo}.  Soon thereafter, B. Kostant proved a version of the Borel-Weil-Bott theorem in which the computation of sheaf cohomology was reduced to Lie algebra cohomology, \cite{Ko}.  The work of R. Langlands and W. Schmid extended some of the results of Borel-Weil, Bott, and Kostant to certain infinite dimensional Harish-Chandra modules, \cite{L1}, \cite{Sc}.  In a further step, cohomological induction emerged as an infinitesimal counter-part to sheaf cohomology, and led to the construction of a broader class of Harish-Chandra modules, [V2], [Z], [EW].  In contrast to the Borel-Weil-Bott theorem, where every simple finite dimensional $\frak g$-module appears as a higher cohomology group, not every simple Harish-Chandra module appears as a submodule of a module cohomologically induced from a proper compatible parabolic subalgebra (the definition of a compatible parabolic subalgebra see in section 1 below).  This observation applies also to generalized Harish-Chandra modules, and therefore the study of the fundamental series is only a first step towards a classification of simple generalized Harish-Chandra modules.  Moreover, the construction and characterization of the fundamental series is merely a branch of a tree whose trunk is the classical Borel-Weil-Bott theorem.

One more common feature of this present work with the work of Armand Borel is that we study general (non-symmetric) reductive pairs $(\frak g,\frak k)$, which have appeared in Borel's work on the topology of homogeneous spaces, \cite{B}.

The paper is organized as follows. In section 0 we fix the notation.  In section 1 we present a minimum of background material and then state the main new results. Theorems 1-4 and Corollaries 1-2.  Sections 2 and 3 are the technical core of the paper; in section 2 we prove Theorem 1, and in section 3 we prove all other results of section 1.  Unfortunately, the proofs are not self-contained as our work relies heavily on the machinery developed by D. Vogan in the course of his fundamental work \cite{V2}.  We state all results in the generality we need them but we often refer to \cite{V2} for the proof if it does not require essential modifications.  Finally in section 4, we discuss some particular cases in our construction, and in particular consider in more detail the case when $\k$ is a $s\ell (2)$-subalgebra.  In this case, the genericity condition on a $\k$-type reduces to a simple explicit inequality.

\noindent{\bf Acknowledgement.}  We thank the referee for pointing out certain inaccuracies.  The first named author gratefully acknowledges support from the NSF, the Max Planck Institute for Mathematics in Bonn and Yale University.

\head 0. Conventions  \endhead

The ground field is $\C$, and if not explicitly stated otherwise, all vector spaces and Lie algebras are defined over $\C$.  By definition, $\N = \{0,1,2,\ldots\}$.  The sign $\otimes$ denotes tensor product over $\C$.  The superscript $*$ indicates dual space.  The sign $\cplus$ stands for semidirect sum of Lie algebras (if $\frak l = \frak l'\cplus\frak l''$, $\frak l'$ is an ideal in $\frak l$ and $\frak l''\simeq\frak l/\frak l'$). $H^\cdot (\frak l, M)$ stands for the cohomology of a Lie algebra $\frak l$ with coefficients in an $\frak l$-module $M$, and $M^\frak l = H^0(\frak l,M)$ stands for space of $\frak l$-invariants of $M$.  $\Lambda^\cdot(\quad)$ and $S^\cdot (\quad)$ denote respectively the exterior and symmetric algebra.

If $\frak l$ is a Lie algebra, then $U(\frak l)$ stands for the enveloping algebra of $\frak l$ and $Z_{U(\frak l)}$ denotes the center of $U(\frak l)$.  We identify $\frak l$-modules with $U(\frak l)$-modules.  It is well known that if $\frak l$ is finite dimensional and $M$ is a simple $\frak l$-module (or equivalently a simple $U(\frak l)$-module), $Z_{U(\frak l)}$ acts on $M$ via a $Z_{U(\frak l)}$-{\it character}, i.e. via an algebra homomorphism $\theta_M :Z_{U(\frak l)} \to\C$.  By $Z(\frak l)$ we denote the center of the Lie algebra $\frak l$.

If $\frak l$ is a Lie algebra, $M$ is an $\frak l$-module, and $\omega\in\frak l^*$, we put \newline $M^\omega := \{m\in M\vert\frak l\cdot m = \omega (\frak l)m \,\,  \forall\ell\in\frak l\}$.  We call $M^\omega$ a {\it weight space} of $M$ and we say that $M$ is an $\frak l$-{\it weight module} if
$$
M = \bigoplus\limits_{\omega\in\frak l^{*}} M^\omega.
$$
By $\supp_{\frak l} M$ we denote the set $\{\omega\in\frak l^*\vert M^\omega\ne 0\}$.

A finite {\it multiset} is a function $f$ from a finite set $D$ into $\N$.  A {\it submultiset} of $f$ is a multiset $f'$ defined on the same domain $D$ such that $f'(d)\leq f(d)$ for any $d\in D$.  For any finite multiset $f$, defined on an additive monoid $D$, we can put $\rho_f := \frac{1}{2}\sum\limits_{d\in D} f(d)d$.  If $M$ is an $\frak l$-weight module as above, and $\dim M < \infty$, $M$ determines the finite multiset  $\ch_{\frak l} M$ which is the function $\omega\mapsto \dim M^\omega$ defined on $\supp_{\frak l} M$.

\head 1.  Statement of results \endhead

\noindent{\bf 1.1. Reductive pairs, compatible parabolics and generic $\frak k$-types.}  Let $\frak g$ be a finite dimensional semisimple Lie algebra and $\frak k\subset\frak g$ be an algebraic subalgebra which is reductive in $\frak g$.  We fix a Cartan subalgebra $\frak t$ of $\frak k$ and a Cartan subalgebra $\frak h$ of $\frak g$ such that $\frak t\subset\frak h$. (If $(\frak g,\frak k)$ is a symmetric pair, then $\frak h$ is unique and is called a fundamental Cartan subalgebra.  An important feature of the general case we consider is that $\frak h$ is no longer unique).  By $\Delta$ we denote the set of $\frak h$-roots of $\frak g$, i.e. $\Delta = \{\supp_{\frak h} \frak g\}\backslash \{ 0\}$.  Note that, since $\frak k$ is reductive in $\frak g$, $\frak g$ is a $\frak t$-weight module.  Therefore we can set $\Delta_{\frak t} := \{\supp_{\frak t}\frak g\} \backslash \{ 0\}$.  Note also that the $\R$-span of the roots of $\frak g$ fixes a real structure on $\frak h^*$, whose projection onto $\frak t^*$ is a well-defined real structure on $\frak t^*$.  In what follows, we will denote by $\Re \lambda$ the real part of an element $\lambda\in\frak t^*$.  We fix also a Borel subalgebra $\frak b_{\frak k} \subset\frak k$ with $\frak b_{\frak k}\supset\frak t$.  Then $\frak b_{\frak k} = \frak t ~\crplus \frak n_{\frak k}$, where $\frak n_{\frak k}$ is the nilradical of $\frak b_{\frak k}$.   We set $\rho :=\rho_{\ch_{\frak t} \frak n_{\frak k}}$.  The quintet $\frak g,\frak h,\frak k,\frak b_{\frak k},\frak t$ will be fixed throughout the paper.  By $W_{\frak k}$ we denote the Weyl group of $\frak k$.

As usual, we will parametrize the characters of $Z_{U(\frak g)}$ via the Harish-Chandra homomorphism.  More precisely, if $\frak b_{\frak g}$ is a given Borel subalgebra of $\frak g$ with $\frak b_{\frak g} \supset\frak h$ ($\frak b_{\frak g}$ will be specified below), the $Z_{U(\frak g)}$-character corresponding to $\kappa\in\frak h^*$ via the Harish-Chandra homomorphism defined by $\frak b_{\frak g}$ will be denoted by $\theta_{\kappa}$ ($\theta_{\rho_{\ch_{\h}\b_{\g}}}$ is then the trivial $Z_{U(\frak g)}$-character).

By $ \langle \, ,\, \rangle$ we denote the unique $\frak g$-invariant symmetric bilinear form on $\frak g^*$ such that $\langle \alpha,\alpha \rangle =2$ for any long root of a simple component of $\frak g$.  The form $\langle \, , \, \rangle$ enables us to identify $\frak g$ with $\frak g^*$. Then $\frak h$ is identified with $\frak h^*$, and $\frak k$ is identified with $\frak k^*$.  We will sometimes consider $\langle \, ,\, \rangle$ as a form on $\frak g$.  The superscript $\perp$ indicates orthogonal space.  Note that there is a canonical $\k$-module decomposition $\frak g = \frak k\oplus\frak k^\perp$.  We also set $\Vert\kappa\Vert^2 := \langle \kappa,\kappa \rangle$ for any $\kappa\in\frak h^*$.

We say that an element $\lambda\in\frak t^*$ is $(\frak g,\frak k)$-{\it regular} if $\langle \Re\lambda,\sigma \rangle\ne 0$ for all $\sigma\in\Delta_{\frak t}$.  To any $\lambda\in\frak t^*$ we associate the following parabolic subalgebra $\frak p_\lambda$ of $\frak g$:
$$
\frak p_\lambda = \frak  h \oplus (\bigoplus\limits_{\alpha\in\Delta_{\lambda}} \frak g^\alpha),
$$
where $\Delta_\lambda :=\{\alpha\in\Delta\mid \langle \Re \lambda, \alpha \rangle\geq 0\}$.  By $\frak m_\lambda$ and $\frak n_\lambda$ we denote respectively the reductive part of $\frak p_\lambda$ (containing $\frak h$) and the nilradical of $\frak p_\lambda$.  In particular $\frak p_\lambda = \m_\lambda ~\crplus \n_\lambda$, and if $\lambda$ is $\frak b_{\frak k}$-dominant, then $\frak p_\lambda\cap \frak k = \frak n_{\frak  k}$.    We call $\frak p_\lambda$ a {\it parabolic subalgebra compatible with} $\frak t$, or simply a {\it compatible parabolic subalgebra}.  A compatible parabolic subalgebra $\frak p = \frak m ~\crplus \frak n$ (i.e. $\frak p = \frak p_\lambda$ for some $\lambda\in\frak t^*$) is {\it minimal} if it does not properly contain another compatible parabolic subalgebra.  It is an important observation that if $\frak p = \frak m ~\crplus\frak n$ is minimal, then $\frak t\subset Z(\frak m)$.  Furthermore, it is easy to see that a compatible parabolic subalgebra $\frak p_\lambda$ is minimal if and only if $\frak m_\lambda$ equals the centralizer $C(\frak t)$ of $\frak t$ in $\frak g$, or equivalently if and only if $\lambda$ is $(\frak g,\frak k)$-regular.

A {\it $\k$-type} is by definition a simple finite dimensional $\k$-module.  By $V(\mu)$ we will denote a $\frak k$-type with $\frak b_{\frak k}$-highest weight $\mu$ ($\mu$ is then $\frak k$-integral and $\frak b_{\frak k}$-dominant).  Let $V(\mu)$ be a $\k$-type such that $\mu + 2\rho$ is  $(\g,\k)$ regular, and let $\p = \m ~\crplus\n$ be the minimal compatible parabolic subalgebra $\p_{\mu +2\rho}$.  Put $\rho_\n := \rho_{\ch_{\t}\n}$. We define $V(\mu)$ to be {\it generic} if the following two conditions hold:

(1)  $\langle \Re\mu + 2\rho - \rho_\n,\alpha \rangle\geq 0 ~ \forall\alpha\in \supp_\t \n_\k$:

(2)  $\langle \Re \mu + 2\rho - \rho_S, \rho_S \rangle > 0$ for every submultiset $S$ of $\ch_\t \n$.
\vskip .15in

\noindent{\bf 1.2. $(\g,\k)$-modules of finite type and minimal $\k$-types.}  For the purposes of this paper we will call a $\g$-module $M$ a $(\g,\k)$-{\it module} if $M$ is isomorphic as a $\k$-module to a direct sum of isotypic components of $\k$-types.  If $M$ is a $(\g,\k)$-module, we write $M[\mu]$ for the $V(\mu)$-isotypic component of $M$, and we say that $V(\mu)$ is a $\k$-{\it type of} $M$ if $M[\mu]\ne 0$.  We say that a $(\g,\k)$-module is {\it of finite type} if $\dim M[\mu] \neq \infty$ for every $\k$-type $V(\mu)$.  We will also refer to $(\g,\k)$-modules of finite type as {\it generalized Harish-Chandra modules}.

Let $\Theta_\k$ be the discrete subgroup of $Z(\k)^*$ generated by $\supp_{Z(\k)} \g$.  By $\Cal M$ we denote the class of $(\g,\k)$-modules $M$ for which there exists a finite subset $S\subset Z(\k)^*$ such that $\supp_{Z(\k)} M\subset S + \Theta_\k$.  If $M$ is a module in $\Cal M$, a $\k$-type $V(\mu)$ of $M$ is {\it minimal} if the function $\mu'\mapsto \Vert\Re \mu' + 2\rho\Vert^2$ defined on the set $\{ \mu'\in\t^*\mid M[\mu']\ne 0\}$ has a minimum at $\mu$.  Any non-zero $(\g,\k)$-module $M$ in $\Cal M$ has a minimal $\k$-type.  This follows from the fact that the squared length of a vector has a minimum on every shifted lattice in Euclidean space.
\vskip .15in

\noindent{\bf 1.3.  Existence of $\n$-cohomology.}  Our first result in this paper is the following analog of a theorem of Vogan, \cite{V1}, \cite{V2}.

\proclaim{Theorem 1}  Let $M$ be a module in $\Cal M$ which has a generic minimal $\k$-type $V(\mu)$.  There is a vector space isomorphism
$$
(M^{\n\cap\k})^\mu\otimes \Lambda^r (\n\cap\k^\perp)^* \cong H^r (\n,M)^{\mu -2\rho^\perp_\n}, \tag3
$$
where $\n := \n_{\mu +2\rho}, \rho^\perp_\n := \rho_{\ch_{\t} (\n\cap\k^\perp)}$, and $r:= \dim(\n\cap\k^\perp)$.  Moreover,
$$
H^i(\n, M)^{\mu-2\rho^\perp_\n} = 0
$$
for $i\ne r$.
\endproclaim
\vskip .15in

\noindent{\bf 1.4. The fundamental series of generalized Harish-Chandra modules.}  Our second result is the following construction of a new series of $(\g,\k)$-modules of finite type which we call the {\it fundamental series of generalized Harish-Chandra modules.}  Recall that the functor of $\k$-locally finite vectors $\Gamma_{\k,\t}$ is a well-defined left exact functor on the category of $(\g,\t)$-modules with values in $(\g,\k)$-modules,
$$
\Gamma_{\k,\t} (M) = \sum\limits_{M'\subset M,\dim M'=1, \dim U(\k)\cdot M' < \infty} M'.
$$
By $R^\cdot \Gamma_{\k,\t}:=\bigoplus\limits_{i\geq 0} R^i \Gamma_{\k,\t}$ we denote as usual the total right derived functor of $\Gamma_{\k,\t}$, see \cite{PZ} and the references therein.

We need also the following ``production'' or ``coinduction'' functor from the category of $(\p,\t)$-modules to the category of $(\g,\t)$-modules:
$$
\pro^{\g,\t}_{\p,\t} (N) := \Gamma_{\t,0} (\Hom_{U(\p)} (U(\g), N)).
$$
The functor $\pro^{\g,\t}_{\p,\t}$ is exact.

We are now ready to state our second theorem, which constructs and describes the {\it fundamental series} of $(\g,\k)$-modules of finite type $F^\cdot (\p,E)$.

\proclaim{Theorem 2} Let $\p = \m ~\crplus\n$ be a minimal compatible parabolic subalgebra, $E$ be a simple finite dimensional $\p$-module on which $\t$ acts via the weight $\omega\in\t^*$, $\rho_\n := \rho_{\ch_{\t}\n},\rho^\perp_{\n} := \rho_{\ch_{\t}(\n\cap\k^{\perp})}$, and $\mu := \omega + 2\rho^\perp_\n$.  Set
$$
F^\cdot(\p, E) := R^\cdot\Gamma_{\k,\t} (\pro^{\g,\t}_{\p,\t} (E\otimes\Lambda^{\dim \n}(\n))).
$$
Then the following assertions hold under the assumptions that $\p = \p_{\mu +2\rho}$ and that $\mu$ is $\b_\k$-dominant, $\k$-integral and yields a generic $\k$-type $V(\mu)$.

a) $F^\cdot(\p,E)$ is a  $(\g,\k)$-module of finite type in the class $\Cal M$, and $Z_{U(\g)}$ acts on $F^\cdot (\p ,E)$ via the $Z_{U(\g)}$-character $\theta_{\nu +\tilde\rho}$ where $\tilde\rho := \rho_{\ch_{\h}\b}$ for some fixed Borel subalgebra $\b$ of $\g$ with $\b\supset\h, \b\subset\p$ and $\b\cap\k = \b_\k$, and where $\nu$ is the $\b$-highest weight of $E$ (note that $\nu\mid_\t =\omega$).

b) $F^i (\p,E) = 0$ for $i\ne s :=\dim\n_\k$ (equivalently $s = \dim (\n\cap\k)$).

c)  There is a $\k$-module isomorphism
$$
F^s (\p,E)[\mu ] \cong \C^{\dim E}\otimes V(\mu),
$$
and $V(\mu)$ is the unique minimal $\k$-type of $F^s (\p,E)$.

d)  Let $\bar F^s (\p,E)$ be the $\g$-submodule of $F^s (\p,E)$ generated by $F^s(\p,E)[\mu]$.  Then $\bar F^s (\p,E)$ is the unique simple submodule of $F^s(\p,E)$, and moreover, $\bar F^s(\p,E)$ is a submodule of any $\g$-submodule of $F^s(\p,E)$.

e)  For any non-zero $\g$-submodule $M$ of $F^s(\p,E)$ there is an isomorphism of $\m$-modules
$$
H^r (\n, M)^\omega \cong E.
$$
\endproclaim

\noindent{\bf 1.5. Reconstruction theorems.}  The results in this subsection constitute our main results.

\proclaim{Theorem 3 (First reconstruction theorem)}  Let $M$ be a simple $(\g,\k)$-module of finite type with a minimal $\k$-type $V(\mu)$ which is generic.  Then $\p := \p_{\mu +2\rho} = \m ~\crplus\n$ is a minimal compatible parabolic subalgebra.  Let $\omega := \mu -2\rho^\perp_\n$ and $E$ be the $\p$-module $H^r(\n, M)^\omega$ with trivial $\n$-action, where $r = \dim (\n\cap\k^\perp)$.  Then $E$ is a simple $\p$-module, the pair $(\p,E)$ satisfies the hypotheses of Theorem 2, and $M$ is canonically isomorphic to $\bar F^s (\p,E)$ for $s= \dim (\n\cap\k)$.
\endproclaim

\proclaim{Corollary 1 (Generic version of a theorem of Harish-Chandra)}  There exist at most finitely many simple $(\g,\k)$-modules $M$ of finite type with a fixed $Z_{U(\g)}$-character such that a minimal $\k$-type of $M$ is generic. (Moreover, each such $M$ has a unique minimal $\k$-type by Theorem 2, c).)
\endproclaim

{\bf Proof of Corollary 1.}  By Theorems 2 and 3, if $M$ is a simple $(\g,\k)$-module of finite type with generic minimal $\k$-type $V(\mu )$ for some $\mu$, then the $Z_{U(\g)}$-character of $M$ is $\theta_{\nu +\tilde\rho}$.  There are finitely many Borel subalgebras $\b$ as in Theorem 2, a); hence, if $\theta_{\nu +\tilde\rho}$ is fixed, there are finitely many possibilities for the weight $\nu$ (as $\theta_{\nu +\tilde\rho}$ determines $\nu +\tilde\rho$ up to a finite choice).  Therefore, there are a finitely many possibilities for the $\p$-module $E$, and hence there are finitely many possibilities for $M$.  \qed

\proclaim{Theorem 4 (Second reconstruction theorem)}  Assume that the pair $(\g,\k)$ is regular, i.e. $\t$ \,\, contains a regular element of $\g$.  Let $M$ be a simple  $(\g,\k)$-module (a priori of infinite type) with a minimal $\k$-type $V(\mu)$ which is generic.  Then $M$ has finite type, and hence by Theorem 3, $M$ is canonically isomorphic to $\bar F^s(\p, E)$ (where $\p,E$ and $s$ are as in Theorem 3).
\endproclaim

\proclaim{Corollary 2}  Let the pair $(\g,\k)$ be regular.

a)  There exist at most finitely many simple $(\g,\k)$-modules $M$ with a fixed $Z_{U(\g)}$-character, such that a minimal $\k$-type of $M$ is generic. All such $M$ are of finite type, (and have a unqiue minimal $\k$-type by Theorem 2, c)).

b)  (Generic version of Harish-Chandra's admissibility theorem)  Every simple $(\g,\k)$-module with a generic minimal $\k$-type has finite type.
\endproclaim

{\bf Proof of Corollary 2}  The proof of a) is as the proof of Corollary 1 but uses Theorem 4 instead of Theorem 3, and b) is a direct consequence of Theorem 4.  \qed

The proofs of Theorem 1-4 depend heavily on adaptations of certain important results of D. Vogan \cite{V2}, from the case of a symmetric pair to the case of a general reductive pair $(\g,\k)$, and are presented in the following sections 2 and 3.

\head 2.  $\n$-cohomology \endhead

In this section we present the minimum material on $\n$-cohomology necessary to outline the proof of Theorem 1.

Let $\p = \m ~\crplus\n$ be a compatible parabolic subalgebra corresponding to an element $\lambda\in\t^*$ (i.e. $\p =\p_\lambda$) which we assume $\k$-regular, and let $M$ be a $(\g,\k)$-module in $\Cal M$.

\proclaim{Proposition 1}  In the category of $\t$-weight modules, there exists a bounded (not necessarily first quadrant) cohomology spectral sequence which converges to $H^\cdot (\n,M)$, with
$$
E^{a,b}_1 = H^{a+b-R(a)} (\n\cap\k, M)\otimes V^*_a,
$$
where $a$ runs over $\{ 0,\ldots ,n\}$ for some $n$, $R$ is a monotonic function on $\{ 0,\ldots, n\}$ with values in $\N$ such that $R(a)\leq a$ and $R(n)=r$, $V_a$ is a $\t$-submodule of $\Lambda^{R(a)} (\n\cap\k^\perp)$ for every $a$, and $V_n = \Lambda^r (\n\cap \k^\perp)$.
\endproclaim

The spectral sequence whose existence is claimed in Proposition 1 is a version of the Hochschild-Serre spectral sequence and is constructed explicitly by Vogan in \cite{V2, Theorem 5.2.2} under the assumption that the pair $(\g,\k)$ is symmetric.  However, as this assumption is not used in the construction, we refer the reader directly to \cite{V2}.

Proposition 1 has the following corollary.

\proclaim{Corollary 3}

a) If $M$ is a $(\g,\k)$-module of finite type, then $H^\cdot (\n,M)$ is an $(\m,\t)$-module of finite type.  Moreover, if $M$ is $Z_{U(\g)}$-finite (i.e. the action of $Z_{U(\g)}$ on $M$ factors through a finite dimensional quotient of $Z_{U(\g)}$) then $H^\cdot (\n,M)$ is $Z_{U(\m)}$-finite.

b) If $\p$ is a minimal compatible parabolic subalgebra and $M$ is a $(\g,\k)$-module of finite type which is in addition $Z_{U(\g)}$-finite, then $H^\cdot (\n,M)$ is finite dimensional.
\endproclaim

\demo{Proof}  a) is a straightforward generalization of \cite{V2, Corollary 5.2.4}.  Part b) follows from a) and from the observation that $\t\subset Z(\m)$ whenever $\p$ is minimal.  Indeed, as $H^\cdot (\n,M)$ is an $(\m,\t)$-module of finite type and $\t\subset Z(\m)$, $H^\cdot (\n,M)$ considered as an $\m$-module is a direct sum of finite dimensional isotypic components.  The fact that $H^\cdot (\n,M)$ is $Z_{U(\m)}$-finite implies that there are only finitely many such components, i.e. that $\dim H^\cdot (\n,M) < \infty$. \qed
\enddemo

\proclaim{Corollary 4}  For each $\kappa\in\t^*$ we have a spectral sequence of vector spaces which converges to $H^\cdot (\n,M)^\kappa$ and whose $E_1$-term is $(E^{a,b}_1)^\kappa$, where $E^{a,b}_1$ is as in Proposition 1.  Moreover, there are (edge) homomorphisms
$$
\pi_i^{\kappa +2\rho^{\perp}_{\n}} : H^i(\n\cap\k, M)^{\kappa +2\rho^\perp_\n}\otimes \Lambda^{r} (\n\cap\k^\perp)^* \to H^{i+r} (\n,M)^\kappa, \tag4
$$
where $i = n + b-r$.
\endproclaim

\demo{Proof}  The fact that the spectral sequence of Proposition 1 is a spectral sequence in the category of weight $\t$-modules implies that it has a well-defined direct summand consisting of $\kappa$-weight vectors.  Its corresponding $E^{a,b}_1$-term equals $(E_1^{a,b})^\kappa$.
\enddemo

In \cite{V2, 5.2} Vogan constructs (under the assumption that the pair $(\g,\k)$ is symmetric) linear maps $(E_1^{n,b})^\kappa \to (E_\infty^{n,b})^\kappa$ which in turn yield edge homomorphisms for the spectral sequence with term $(E_1^{a,b})^\kappa$,
$$
\pi_i^{\kappa +2\rho^\perp_\n} : H^i (\n\cap\k, M)^{\kappa +2\rho^\perp_\n} \otimes\Lambda^r (\n\cap\k^\perp)^* \to H^{i+r} (\n,M)^\kappa \cong\bigoplus\limits_{a+b=i+r} (E^{a,b}_\infty )^\kappa.
$$
This construction extends to the more general case we consider, and we refer the reader to \cite{V2}. \qed

\proclaim{Proposition 2}  Under the hypothesis of Theorem 1, $(E^{a,b}_1)^{\mu -2\rho^\perp_\n} = 0$ for $(a,b)\ne (n,r-n)$; therefore the spectral sequence from Corollary 4 for $\kappa = \mu -2\rho^\perp_\n$ collapses at the level $E_1$. 
\endproclaim

\demo{Proof}  Let $V(\delta )$ be a $\k$-type such that
$$
H^i (\n\cap\k, V(\delta ))\otimes \Lambda^j (\n\cap\k^\perp)^* \ne 0
$$
for some $(i,j)\ne (0,r)$.  Then, according to \cite{V2, 5.4} there exist $\sigma\in W_\k$ and a submultiset $A$ of $\ch_\t (\n\cap\k^\perp)$ such that $\sigma (\delta +\rho ) -\rho = \mu - 2\rho_A$.  Hence $\sigma (\delta +2\rho) = \mu + \rho -2\rho_A + \sigma\rho = \mu + 2\rho - (2\rho_A +\rho -\sigma\rho ) = \mu + 2\rho - 2\rho_B$ for an appropriate submultiset $B$ of $\ch_\t \n$.  Furthermore, as $V(\mu)$ is generic, $\langle \Re \mu + 2\rho - \rho_B, \rho_B\rangle > 0$ by condition (2), and thus
$$
\Vert\Re\delta + 2\rho\Vert^2 = \Vert\sigma (\Re\delta + 2\rho)\Vert^2 = \Vert\Re\mu + 2\rho\Vert^2
- 4 \langle \Re\mu +2\rho - \rho_B, \rho_B \rangle  < \Vert\Re\mu + 2\rho\Vert^2.
$$
As $V(\mu)$ is minimal in $M$, $V(\delta )$ is not a $\k$-type of $M$, and hence
$$
(H^i (\n\cap\k. M)\otimes \Lambda^j (\n\cap\k^\perp)^*)^{\mu -2\rho^\perp_\n} = 0
$$
for $(i,j)\ne (0,r)$.  Since $(E_1^{a,b})^{\mu -2\rho^\perp_\n}\subset (H^i (\n\cap\k ,M)\otimes \Lambda^j (\n\cap\k^\perp)^*)^{\mu -2\rho^\perp_\n}$ for \newline $i = a+b-R(a)$ and $j = R(a)$, we obtain
$$
(E_1^{a,b})^{\mu -2\rho^\perp_\n} = 0
$$
for $(a,b)\ne (n, r-n)$. \qed
\enddemo

Theorem 1 follows from Proposition 2 via the observations:  first,
$$
H^i(\n,M)^{\mu -2\rho^\perp_\n} \cong \bigoplus\limits_{a+b=i} (E_\infty^{a,b})^{\mu -2\rho^\perp_\n} = \bigoplus\limits_{a+b=i} (E_1^{a,b})^{\mu -2\rho^\perp_\n}
$$
and hence $H^i (\n,M)^{\mu -2\rho^\perp_\n} = 0$ for $i\ne r$, and second, the desired isomorphism (3) is nothing but the edge isomorphism $\pi^\mu_0$. \qed

\head 3.  Construction and characterization of the fundamental series \endhead

In this section we prove Theorems 2-4.  We have split the proof into several statements, some of which are of interest by themselves.  Throughout this section, $\p, E, \rho^\perp_\n, \mu,\omega$ and $s$ are as in Theorem 2 (with the hypotheses of Theorem 2 holding) and are assumed fixed.

\proclaim{Proposition 3}    Let $V(\delta)$ be a $\k$-type of $F^{s-i}(\p,E)$ for some $i\in\Z$.

a)  There exists $w\in W_\k$ of length $i$ (in particular, $i\in\N$) and a multiset
$$\align
&n_\cdot := \supp_\t (\n\cap\k^\perp) \to \N,\\
&\qquad \beta\mapsto n_\beta
\endalign
$$
such that
$$
\omega = w(\delta +\rho )-\rho -2\rho^\perp_\n - \sum\limits_{\beta} n_\beta \beta. \tag5
$$
Furthermore, the multiplicity of $V(\delta)$ in $F^{s-i} (\p,E)$ is bounded by the integer
$$
\dim E \sum\limits_{\ell (w)=i} \dim (S^\cdot(\frak n\cap \frak k^\perp)^{\xi (w)}),
$$
where $\xi (w)$ is the $\frak t$-weight $w(\delta +\rho)-\rho -\omega -2\rho^\perp_\n$ in $S^\cdot (\n\cap\k^\perp)$.

b) The equality
$$
\sum\limits_{0\leq i\leq s} (-1)^i \dim \Hom_\k (V(\delta), F^{s-i} (\p, E))\tag6 
$$
$$
= \sum\limits_{0\leq j\leq s} (-1)^j (\sum\limits^\infty_{m=0} \dim \Hom_\t (H^j (\n\cap\k, V(\delta)), S^m(\n\cap\k^\perp)\otimes E\otimes\Lambda^{\dim (\n\cap\k^{\perp})} (\n\cap\k^\perp)))
$$
holds, and the inner sum on the right hand side of (6) is finite.
\endproclaim

Proposition 3 is a modification of Vogan's Theorem 6.3.12 and Corollary 6.3.13 in \cite{V2}, and its proof follows exactly the same lines (an inspection of Vogan's proofs reveals that the symmetry assumption on $(\g,\k)$ is not needed).  Therefore, we refer the reader to \cite{V2}.

Proposition 3 implies that $F^\cdot (\p,E)$ is a $(\g,\k)$-module of finite type, and also that $F^i(\p,E) = 0$ for $i > s$.  Furthermore, Proposition 3 implies that $F^\cdot (\p,E)$ is in the class $\Cal M$.  To see this, one has to restrict (5) to $Z(\k)$ and notice that the one element set $\{\omega_{\mid Z(\k)}\}$ can be chosen as the finite set $S$ that appears in the definition of the class $\Cal M$.

\proclaim{Proposition 4}  Suppose $V(\delta)$ is a $\k$-type of $F^{s-i}(\p, E)$ and
$$
(H^\cdot (\n\cap\k, V(\delta))\otimes\Lambda^\cdot (\n\cap\k^\perp)^*)^\omega \ne 0. \tag7
$$
Then $i=0$.
\endproclaim

\demo{Proof}  By Proposition 3, a) there exist $w$ and $n_{\cdot}$ such that (5) holds.  Since $\delta$ is dominant with respect to $\b_\k$,
$$
w\delta = \delta - \sum\limits_{\alpha\in\supp_{\t}(\n\cap\k)} m_\alpha \alpha
$$
for some $m_\alpha\in\R, m_\alpha\geq 0$.  Thus we can rewrite (5) as
$$
\delta = \mu +\rho -w\rho + \sum\limits_{\alpha\in\supp_{\t}(\n\cap\k)} m_\alpha \alpha + \sum\limits_{\beta\in\supp_{\t}(\n\cap\k^{\perp})} n_\beta \beta . \tag8
$$
\enddemo
On the other hand, assumption (7) implies via Kostant's theorem, \cite{Ko}, the existence of an element $\sigma\in W_\k$ and a submultiset $q_\cdot$ of $\ch_\t (\n\cap\k^\perp)$ such that
$$
\sigma (\delta +\rho )-\rho - 2\rho^\perp_\n + \sum\limits_{\gamma\in\supp_\t (\n\cap\k^{\perp})} q_\gamma\gamma = \omega,
$$
or equivalently
$$
\sigma (\delta +\rho ) = \omega + \rho_\n + (\rho^\perp_\n - \sum\limits_{\gamma\in\supp_\t (\n\cap\k^{\perp})} q_\gamma \gamma ), \tag9
$$
as $\rho_\n = \rho + \rho^\perp_\n$.  Notice that
$$
\rho^\perp_\n - \sum\limits_{\gamma\in\supp_\t (\n\cap\k^{\perp})} q_\gamma \gamma = \sigma^{-1} (\rho^\perp_\n -\sum\limits_{\gamma\in\supp_\t (\n\cap\k^{\perp})} q'_\gamma \gamma )
$$
for an appropriate submultiset $q'_\cdot$ of $\supp_\t (\n\cap\k^\perp)$.  Moreover, the genericity condition (1) on $\mu$, rewritten in terms of $\omega$, reads $\langle \Re\omega + \rho_\n, \alpha \rangle \geq 0$ for all $\alpha\in \supp_\t (\n\cap\k)$.  Hence,
$$
\sigma^{-1} (\Re\omega + \rho_\n) = \Re\omega + \rho_\n - \sum\limits_{\alpha\in\supp_\t(\n\cap\k)} p_\alpha \alpha
$$
for some $p_\alpha \in\R$, $p_\alpha\geq 0$.  Since $\omega + 2\rho^\perp_\n$ is $\b_\k$-dominant, $\langle \Im\omega, \alpha \rangle = 0$ for  $\alpha\in\supp_\t (\n\cap\k)$, and hence $\sigma^{-1} (\Im ~\omega) = \Im ~\omega$.  All of this allows us to rewrite (9) as
$$
\delta  =\omega + \rho_\n -\rho - \sum\limits_{\alpha\in\supp_\t(\n\cap\k)} p_\alpha \alpha + \rho^\perp_\n - \sum\limits_{\gamma\in\supp_\t (\n\cap\k^\perp)} q_\gamma \gamma .\tag10 
$$

By comparing (8) and (10), we obtain
$$\align
&\rho -w\rho + \sum\limits_{\alpha\in\supp_\t (\n\cap\k)} m_\alpha \alpha + \sum\limits_{\beta\in\supp_\t (\n\cap\k^{\perp})} n_\beta \beta \tag 11 \\
&= - \sum\limits_{\alpha\in\supp_\t (\n\cap\k)} p_\alpha \alpha -\sum\limits_{\beta\in\supp_\t (\n\cap\k^{\perp})} q'_\gamma \gamma .
\endalign
$$
Since $\langle\Re \mu + 2\rho, \eta\rangle > 0$ for every $\eta\in\supp_\t \n$, (11) implies \newline $\langle\Re\mu + 2\rho, \rho - w\rho \rangle \leq 0$.  As $\rho -w\rho = \sum\limits_{\alpha\in\supp_\t (\n\cap\k)} n_\alpha \alpha$ for $n_\alpha\in\N$, we obtain $\langle \Re\mu + 2\rho, \rho - w\rho \rangle = 0$, or equivalently $\rho = w\rho$.  Therefore, $w = \id$ and, since $i$ is the length of $w$, $i = 0$. \qed

\proclaim{Proposition 5}  (Analog of Frobenius reciprocity)  For any $(\g,\k)$-module $M$ there exist two first quadrant cohomology spectral sequences (in the category of vector spaces) I and II with common limit $\Ext^{a+b}_{\g,\t} (M, \pro^{\g,\t}_{\p,\t} (E\otimes \Lambda^{\dim (\n)} (\n)))$ and respective $E^{a,b}_2$-terms:
$$\align
I^{a,b}_2 &= \Ext^a_{\g,\k} (M, F^b (\p,E)), \\
II^{a,b}_2 &= \Ext^a_{\m,\t} (H^{\dim \n-b} (\n, M), E).
\endalign
$$
\endproclaim

The proof is the same as the proof \cite{V2, Proposition 6.3.2} and we omit it.

We are now able to complete the proof of Theorem 2, b).  We have already shown that $F^i (\p,E) = 0$ for $i>s$  and now we will show that $F^i(\p, E) = 0$ also for $i<s$.  Suppose $F^\cdot (\p,E)\ne 0$ and let $l_0$ be the minimal integer with $F^{l_{0}} (\p,E)\ne 0$.

Set $M = F^{l_{0}} (\p,E)$.  Then
$$
\Hom_{\g,\k} (M, F^{l_{0}}(\p,E))\ne 0,
$$
and $I^{a,b}_2 = 0$ for $a+b < l_0$.  Therefore the spectral sequence $I$ yields an isomorphism $I^{0,l_{0}}_2\to  I^{0,l_{0}}_\infty$, hence $I_\infty^{0,l_{0}}\ne 0$.  As $I$ and $II$ have the same limit, $\bigoplus\limits_{a+b=l_0} I^{a,b}_\infty = \bigoplus\limits_{a+b=l_0} II^{a,b}_\infty$.  Thus $II^{a_0,b_0}_\infty\ne 0$ for some $a_0,b_0$ with $a_0 + b_0=l_0$, and consequently
$$
\Ext^{a_{0}}_{\m,\t} (H^{\dim \n-b_{0}} (\n,M),E)\ne 0. \tag12
$$

By Corollary 3, b) $D:= H^{\dim \n-b_0} (\n,M)$ is a finite dimensional $\t$-weight $\m$-module.  $D$ has a canonical decomposition as
$$
\bigoplus\limits^x_{j=1} C_j\otimes B_j
$$
where $C_j$ are simple non-isomorphic (finite dimensional) modules over the semisimple part $\m_{ss}$ of $\m$ and $B_j$ are (finite dimensional not necessarily semisimple) $Z(\m)$-modules.  Similarly, we can factor $E$ as $C\otimes B$, where $C$ is a simple $\m_{ss}$-module and $B$ is a $1$-dimensional $Z(\m)$-module.  By the Kunneth formula, \cite{We},
$$
\Ext^{a}_{\m,\t} (D,E) =\bigoplus\limits_{1\leq j\leq x, p+q=a} \Ext^p_{\m_{ss}} (C_j,C)\otimes \Ext^q_{Z(\m),\t} (B_j,B).
$$
Furthermore, by Whitehead's lemma, \cite{We},
$$
\Ext^\cdot_{\m_{ss}} (C_j, C) = 0
$$
if $C_j$ and $C$ are inequivalent.  Thus (12) implies that $C\simeq C_{j_{0}}$ for exactly one value $j_0$ of $j$ and
$$
\Ext^{a_{0}}_{\m,\t} (D,E)\simeq \bigoplus\limits_{p+q=a_{0}} \Ext^p_{\m_{ss}} (C_{j_{0}}, C)\otimes \Ext^q_{Z(\m),\t} (B_{j_{0}}, B).
$$
The non-vanishing of $\Ext^\cdot_{Z(\m),\t} (B_{j_{0}}, B)$ implies that the $1$-dimensional $Z(\m)$-module $B$ is a quotient of $B_{j_{0}}$, and hence the $\m$-module $E = C\otimes B$ is a quotient of the $\m$-module $D\simeq C\otimes B_{j_{0}}$.  Therefore we can now conclude that
$$
(H^{\dim \n-b_{0}} (\n,M))^\omega \ne 0 .
$$

As a next step, we apply the spectral sequence from Proposition 1 to show that
$$
(H^\cdot (\n\cap\k, M), V(\delta))\otimes \Lambda^\cdot (\n\cap\k^\perp)^*)^\omega \ne 0,
$$
and we complete the proof of Theorem 2, b) by applying Proposition 4 which yields $s-l_0 = 0$, i.e. $l_0=s$.  \qed

Next we prove assertion c) of Theorem 2.  Theorem 2, b) enables us to rewrite (6) as
$$\align
&\dim \Hom_\k (V(\delta), F^{s}(\p,E)) \\
&=\sum\limits_{0\leq j\leq s} (-1)^j (\sum\limits^\infty_{m=0} \dim \Hom_\t (H^j(\n\cap\k, V(\mu)),
S^m(\n\cap\k^\perp)\otimes E\otimes \Lambda^{\dim (\n\cap\k^{\perp})} (\n\cap\k^\perp))),
\endalign
$$
and, by Kostant's theorem, $\supp_\t H^\cdot (\n\cap\k, V(\mu)) =  \{\tilde\sigma (\mu +\rho )-\rho\mid\tilde \sigma\in W_\k\}$ and $\mu$ appears with multiplicity $1$ in $\{\tilde\sigma (\mu +\rho)-\rho\mid\tilde\sigma\in W_\k\}$.  On the other hand
$$\align
&\supp_\t (S^\cdot (\n\cap\k^\perp)\otimes E\otimes \Lambda^{\dim (\n\cap\k^{\perp})} (\n\cap\k^\perp)) \\
&\qquad  = \{\mu + \sum\limits_{\beta\in\supp_\t (\n\cap\k^{\perp})} n_\beta ~\beta \mid n_\beta\in\N\}.
\endalign
$$
Since $\p=\p_{\mu +2\rho}$, $\langle\Re\mu +2\rho,\alpha \rangle > 0 ~\forall\alpha\in\supp_t\n$; hence
$$
\{\tilde \sigma (\mu +\rho)-\rho\mid\tilde\sigma\in W_\k \}\subset \{ \mu -\sum\limits_{\beta\in\supp_\t (\n\cap\k^\perp)} m_\beta \beta\mid m_\beta\in\N\}.
$$
Therefore,
$$
\{\tilde\sigma (\mu +\rho )-\rho\mid\tilde\sigma\in  W_\k\}\cap\{\mu + \sum\limits_{\beta\in\supp_\t (\n\cap\k^\perp)}\} = \{\mu \},
$$
and consequently
$$ 
\Hom_\t (H^j (\n\cap\k, V(\mu)), S^m (\n\cap\k^\perp)\otimes E\otimes \Lambda^{\dim (\n\cap\k^\perp)} (\n\cap\k^\perp ))\ne 0
$$
only for $m=0$.  This shows that
$$\align
&\dim\Hom_\k (V(\mu), F^s (\p,E))  \\
&= \dim \Hom_\t (H^0 (\n\cap\k, V(\mu)), E\otimes \Lambda^{\dim(\n\cap\k^\perp)} (\n\cap\k^\perp)) = \dim E,
\endalign
$$
i.e. that $V(\mu)$ is a $\k$-type of $F^s (\p,E)$ with multiplicity $\dim E$.

Furthermore, if $F^s(\p,E)[\delta] \ne 0$ for some $V(\delta), \delta\ne\mu$, equality (8) holds with $w=\id$ by Proposition 4, i.e.
$$
\delta = \mu + \sum\limits_{\alpha\in\supp_\t (\n\cap\k)} m_\alpha \alpha + \sum\limits_{\beta\in\supp_\t (\n\cap\k^\perp)}  n_\beta\beta.
$$
Hence $\langle\Re\mu + 2\rho, \alpha\rangle > 0 ~\forall\alpha\in\supp_\t\n$ implies
$$
\Vert\Re\delta + 2\rho\Vert^2 > \Vert\Re \mu + 2\rho\Vert^2,
$$
i.e. $V(\mu )$ is the unique minimal type of $F^s(\p,E)$.  This completes the proof of Theorem 2, c). \qed

\proclaim{Proposition 6}  Let $M$ be a $(\g,\k)$-module.  There exists a (not necessarily first quadrant) cohomology spectral sequence with $E_2$-term
$$
E^{a,b}_2 = \Ext^a_{\m,\t} (H^{r-b} (\n, M),E)
$$
converging to
$$
\Ext^{a+b}_{\g,\k} (M, F^s (\p,E)).
$$
If, in addition, $M$ is a $(\g,\k)$-module of finite type on which $Z_{U(\g)}$ acts via a character, the spectral sequence is a first quadrant spectral sequence (i.e. $E_2^{a,b} = 0$ for $b\leq 0$), and the corner isomorphism $E_2^{0,0}\to E_\infty^{0,0}$ yields an isomorphism
$$
\Hom_\g (M, F^s(\p,E))\simeq \Hom_\m (H^r(\n, M),E). \tag13
$$
\endproclaim

\demo{Proof} The existence is proved by essentially the same argument as in the proof of \cite{V2, Corollary 6.3.4} and uses Proposition 2 and Theorem 2, b).  \qed
\enddemo

If $M$ is a $(\g,\k)$-module of finite type which affords a $Z_{U(\g)}$-chracter, $H^\cdot (\n,M)$ is finite dimensional by Corollary 3, b).  Choose $b_0$ to be the least possible integer with
$$
\Ext^\cdot_{\m,\t} (H^{r-b_0} (\n, M), E)\ne 0.
$$
By the same argument as in the proof of Theorem 2, b), we conclude that
$$
\Hom_\m (H^{r-b_0} (\n,M),E)\ne 0.
$$
Thus $E_2^{0,b_0} \ne 0$ and $E^{a,b}_2 =0$ for $b < b_0$.  Consequently $E_2^{0,b_0}\simeq E^{0,b_0}_\infty$, and we deduce that $\Ext^{b_0}_{\g,\k} (M, F^s(\p,E))\ne 0$.  This enables us to conclude that the spectral sequence is a first quadrant spectral sequence as $b_0\geq 0$, and thus the corner isomorphism yields the desired isomorphism (13).  (Compare our proof with the proof of Theorem 6.5.9, f) in \cite{V2}) \qed

\proclaim{Corollary 5}  If $M$ is a submodule of $F^s(\p,E)$, then $(H^r(\n,M))^\omega\ne 0$.
\endproclaim

\proclaim{Proposition 7}  Let $M$ be a $(\g,\k)$-module with the property that $M[\delta]=0$ for all $\delta$ with $\Vert\Re\delta +2\rho\Vert^2 < \Vert\Re\mu +2\rho\Vert^2$.  Then the isomorphism (3) holds, and in particular $M[\mu]\ne 0$ if and only if $H^r(\n,M)^\omega\ne 0$.
\endproclaim

\demo{Proof} The statement is a consequence of the proof of Theorem 1. \qed
\enddemo

\proclaim{Proposition 8}  For every submodule $M\subset F^s(\p,E), M[\mu] \ne 0$.
\endproclaim

\demo{Proof}  The statement is a direct consequence of Theorem 2, c), Corollary 5 and Proposition 7. \qed
\enddemo

We are now ready to prove Theorem 2, d).  We start with the remark that, if $M$ is any $(\g,\k)$-module of finite type, and $M^*_\k$ is its $\k$-finite dual, i.e. $M^*_\k = \Gamma_{\k,0} (M^*)$, then $M[\mu]^*$ is a $\k$-isotypic component of $M^*_\k$.  Consider the $(\g,\k)$-module of finite type $M := F^s(\p,E)$ and note that Proposition 7 implies that $F^s(\p,E)^*_\k$ is generated by its isotypic component $M[\mu]^*$.  Indeed, if $X$ is the submodule of $F^s(\p,E)^*_\k$ generated by $M[\mu]^*$, and $Y$ is the submodule of $F^s(\p,E)$ orthogonal to $X$, then $Y[\mu] = 0$.  Hence $Y=0$ by Proposition 8, i.e. $X = F^s (\p,E)^*_\k$.

Since $F^s(\p, E)^*_\k$ is generated by $M[\mu]^*, F^s(\p,E)^*$ is of course finitely generated, and as $U(\g)$ is a left Noetherian algebra, $F^s(\p, E)^*_\k$ is a Noetherian $\g$-module.  Therefore $F^s(\p,E)$ is an Artinian $\g$-module. Denote by $\bar F^s(\p,E)$ the $\g$-submodule of $F^s(\p,E)$ generated by $F^s(\p, E)[\mu]$.  Then $\bar F^s(\p,E)$ is both Noetherian and Artinian, and hence by a standard argument in module theory, has finite length.

Let $M_1$ be a simple submodule of $\bar F^s(\p,E)$.  By (13), there is a non-zero $\m$-module map of $H^r(\n, M_1)^\omega$ onto $E$.  By Proposition 7, there is an isomorphism
$$
H^0(\n\cap\k, M_1)^\mu\otimes \Lambda^r(\n\cap\k^\perp)^*\cong H^r(\n,M_1)^\omega ,
$$
and hence $\dim H^0(\n\cap\k, M_1)^\mu\geq\dim E$.  But we also know that
$$
\dim H^0 (\n\cap\k, F^s(\p,E)) = \dim E,
$$
by Theorem 2, c).  Thus, 
$$
\dim H^0 (\n\cap\k, M_1)^\mu = \dim E
$$
and hence
$$
M_1[\mu] = F^s (\p, E)[\mu].
$$
We conclude that
$$
M_1 = \bar F^s(\p, E),
$$
and the proof of Theorem 2, d) is complete. \qed 

The proof of Theorem 2, e) is similar. 

We are now ready to prove Theorem 3.  Since $M$ is a simple $(\g,\k)$-module of finite type, $M$ is in the class $\Cal M$ and $Z_{U(\g)}$ acts on $M$ via a character.  Therefore, by Corollary 3, $H^\cdot (\n,M)$ is a finite dimensional $\t$-weight $\m$-module.  By Theorem 1, $H^r(\n,M)^\omega \ne 0$.  Let $\tilde E$ be any simple quotient of the $\m$-module $H^r(\n,M)^\omega$.  Consider $\tilde E$ as a simple $\p$-module by letting $\n$ act trivially on $\tilde E$.

The fact that $\mu$ is generic implies that the pair $(\p,\tilde E)$ satisfies the hypotheses of Theorem 2.  Thus $F^s (\p,\tilde E)\ne 0$ and there is a canonical isomorphism
$$
\Hom_\g (M, F^s (\p,\tilde E))\cong \Hom_\m (H^r(\n, M),\tilde E). \tag14
$$
Hence, the surjection $H^r(\n, M)\to\tilde E$ determines, via (14), a canonical $\g$-module isomorphism
$$
M\cong \bar F^s (\p,\tilde E). \tag15
$$
Therefore, by Theorem 2, e) $\tilde E$ is isomorphic to $H^r(\n, M)^\omega$, and the surjection \newline $H^r(\n,M)^\omega\to\tilde E$ can be chosen as the identity map.  This implies finally that the isomorphism (15) is a canonical isomorphism
$$
M\cong \bar F^s (\p,E)
$$
for $E = H^r (\n, M)^\omega$, as required.  \qed

Here is the proof of Theorem 4.  The regularity of $\k$ in $\g$ implies that $\m$ equals the Cartan subalgebra $\h$.  Since $M$ is simple, $M$ affords a $Z_{U(\g)}$-character, and the Casselman-Osborne theorem, \cite{CO}, implies that $U(\h) = S^\cdot (\h)$ acts on $H^\cdot (\n,M)$ through a finite dimensional quotient $Q_\cdot$.  Let $J$ be the radical of the algebra $Q$.  Since $J$ is nilpotent, for any non-zero $Q$-module $Z$, $Z/JZ\ne 0$, and hence (by Zorn's Lemma) $Z$ has a $1$-dimensional quotient.  By Theorem 1, $H^r(n,M)^\omega$ is a non-zero $Q$-module (possibly infinite dimensional).  Let $\tilde E$ be a $1$-dimensional quotient of $H^r(\n, M)^\omega$.  The pair $(\p,\tilde E)$ satisfies the hypothesis of Theorem 2 and hence $F^s(\p, \tilde E)$ is non-zero, whereas $F^i(\p,\tilde E) = 0$ for $i\ne s$.

By \cite{V2, Corollary 6.3.4}, we have a bounded spectral sequence with $E_2$-term
$$
\Ext^a_{\h,\t} (H^{(\dim\n-s)-b} (\n, X),\tilde E)
$$
which converges to
$$
\Ext^{a+b}_{\g,\k} (X, F^s(\p, \tilde E))
$$
for any $(\g,\k)$-module $X$.  Set $X = M$ in the above and recall that $\dim \n -s = r$.  Our spectral sequence becomes
$$
\Ext^a_{\h,\t}(H^{r-b}(\n, M)^\omega, \tilde E)\Rightarrow \Ext^{a+b}_{\g,\k} (M, F^s(\p,\tilde E)).
$$

Following Theorem 6.5.9 and its proof in \cite{V2}, choose $b_0$ to be the least possible integer with
$$
\Ext^\cdot_{\h,\t} (H^{r-b_{0}} (n, M)^\omega, \tilde E)\ne 0.
$$
Let $I$ be the maximal ideal in $Q$ which annihilates $\tilde E$.  By elementary homological algebra,
$$
\Ext^\cdot _{\h,\t} ((H^{r-b_{0}} (\n, M)^\omega )_I, \tilde E)\ne 0
$$
where the subscript $I$ indicates localization at $I$.  But then
$$
E^{0,b_0}_2 = \Hom_\h (H^{r-b_{0}} (n, M)^\omega, \tilde E)\ne 0,
$$
and by the assumption on $b_0$, $E^{a,b}_2 = 0$ for any $b < b_0$.  So $E^{0,b_{0}} _2\cong E_\infty^{0,b_{0}}$, and we deduce that
$$
\Ext^{b_{0}}_{\g,\k} (M, F^s (\p, \tilde E))\ne 0.
$$
Thus, $b_0 \geq 0$ and our spectral sequence is a first quadrant spectral sequence.  The corner isomorphism becomes
$$
\Hom_\g (M, F^s(\p, \tilde E))\cong \Hom_\h (H^r(\n, M),\tilde E),
$$
and by the choice of $\tilde E$,  the right hand side is nonzero.  Thus, $M\cong\bar F^s (\p, E)$.  This completes the proof of Theorem 4. \qed

\head 4.  Discussion and examples \endhead

The results of this paper are well known when $(\g,\k)$ is a symmetric pair.  More precisely, in this case, Theorem 4 is the famous Harish-Chandra admissibility theorem and holds without the genericity condition on $V(\mu)$ (in addition, the regularity assumption for $(\g,\k)$ is automatic in this case), and Theorems 1-3 are results of Vogan and are proved in \cite{V2} under less restrictive conditions than the genericity of $V(\mu)$.

In the case when $\k = \h$ is a Cartan subalgebra of $\g$, there exists a classification of simple $(\g,\k)$-modules of finite type, $[M]$, and in principle the results of the present paper can be derived from the classification.  For instance, for $\k = \h$, Theorem 3 claims that a simple weight module of finite type $M$ with generic minimal weight is a $\b$-lowest weight module, where the Borel subalgebra $\b = \p$ is the minimal compatible parabolic subalgebra of Theorem 3.  Consequently, $M$ is a highest weight module with respect to the opposite Borel subalgebra (which contains $\h$).  Theorem 1 then becomes a statement about the $\n$-cohomology of maximal degree $r = \dim \n$, and via Poincar\'e-duality this is equivalent to the obvious statement about the $\n$-covariants of the simple highest weight module.  It seems however, that our genericity condition has not been previously singled out as a sufficient condition for a simple weight module of finite type to be a lowest (or highest) weight module.  Finally, Theorem 4 does not follow from the classification of simple weight modules of finite type, but in principle it could be derived from the classification of all supports of simple weight modules of infinite type given in \cite{DMP}.

The results of this paper are new in all cases when $\k$ is not a symmetric or a Cartan subalgebra of $\g$.  As a simple illustration, we will conclude the paper by a brief discussion of the case when $\k$ is an $s\ell (2)$-subalgebra.

If $\k\simeq s\ell (2)$, then $\dim\t = 1$, and for any $\alpha\in\Delta_\t, \alpha = a\rho$, where $a:=\frac{\langle\alpha,\rho\rangle}{\Vert\rho\Vert^{2}} > 0$. Moreover, $\mu = \omega + 2\rho^\perp_\n = m\rho$ for $m\in\N$, and the genericity condition is equivalent to the single inequality
$$
\langle\mu + 2\rho - \rho_\n, \rho\rangle > 0,
$$
or to the inequality
$$
\dim V(\mu) = m+1\geq \tilde\rho (h),
$$
where $h$ is the semisimple element of the canonical  $s\ell (2)$-basis $e,f,h$ in $\k$ with $h\in\h$ (recall that $\tilde\rho :=\rho_{\ch_{\h}\b}$).  The integer $\tilde\rho (h)$ depends on the pair $(\g,\k)$ and can be computed in the following way.

Write $\tilde\rho = \Sigma r_i\alpha_i$, where $\alpha_i\in\h^*$ are the simple roots of $\b$.  The non-negative half-integers $r_i$ are well-known, see \cite{B}.  Furthermore, a result of E. Dynkin, \cite{D}, states that $\alpha_i(h)\in\{ 0,1,2\}$, and that $\k$ is a principal $s\ell (2)$-subalgebra if $\alpha_i (h) = 2$ for all $i$.  The final inequality equivalent to the genericity of $V(\mu)$ becomes
$$
\dim V(\mu) = m+1\geq \sum\limits_{i} \alpha_i(h)r_i.
$$
In particular, for a principal $s\ell (2)$-subalgebra it reads
$$
m+ 1\geq 2 (\sum\limits_i r_i). \tag16
$$

If $\g = s\ell (3)$ and $\k$ is a principal $s\ell (2)$-subalgebra, the pair $(\g,\k)$ is nothing but the symmetric pair $(s\ell (3), so(3))$ and (16) is the inequality
$$
m\geq 3,
$$
which is well-known to be the necessary and sufficient condition for the first reconstruction theorem to hold.  If $\g = so(5)$ and $\k$ is a principal $s\ell (2)$-subalgebra, (16) is equivalent to
$$
m\geq 6,
$$
and the case when $m\leq 5$ is the ``smallest'' case when the problem of classifying all simple $(\g,\k)$- modules of finite type is still open.

\vskip .20in

\noindent Ivan Penkov

\noindent Department of Mathematics

\noindent University of California-Riverside

\noindent Riverside, California  92521

\noindent email:  penkov\@math.ucr.edu
\vskip .20in

\noindent Gregg Zuckerman

\noindent Department of Mathemaics

\noindent Yale University

\noindent 10 Hillhouse Avenue

\noindent P.O. Box 208283

\noindent New Haven, CT 06520-8283, USA

\noindent email:  gregg\@math.yale.edu

\end
\Refs
\vskip .15in

\ref \key BB \by A. Beilinson, J. Bernstein\pages
\paper Localisation de $\g$-modules
\jour C.R. Acad. Sc. Paris 
{\bf 292} \yr1981
\endref
\vskip .08in

\ref \key B \by A. Borel \pages 115-207
\paper Sur la cohomologie des espaces fibr\'es principaux et des espaces homogen\'enes de groups de Lie compacts
\jour Ann. of Math. (2) {\bf 57} 
\yr1953
\endref
\vskip .06in

\ref \key Bo \by R. Bott \pages 203-248
\paper Homogeneous vector bundles
\jour Ann. of Math.
{\bf 66} \yr 1957
\endref
\vskip .06in

\ref \key Bou \by N. Bourbaki \pages
\paper Lie Groups and Lie Algebras
\paperinfo Ch. 4-6 
\jour Elements of Mathematics, Springer 
\yr 2002
\endref
\vskip .06in

\ref \key CO \by W. Casselman, M.S. Osborne\pages 219-227
\paper The $\n$-cohomology of representations with an infinitesimal character
\jour Comp. Math
{\bf 31} \yr 1975
\endref
\vskip .06in

\ref \key DMP \by I.Dimitrov, O. Mathieu, I. Penkov \pages 2857-2869
\paper On the structure of weight modules
\jour Trans. Amer. Math. Soc.
{\bf  352} \yr 2000
\endref
\vskip .06in

\ref \key D\by E. Dynkin \pages 111-244
\paper Semisimple subalgebras of semisimple Lie algebras
\paperinfo Mat. Sbornik (N.S.)  {\bf 30} (72) (1952), 349-462 (Russian); English: Amer. Math. Soc. Transl. {\bf 6} (1957)
\endref
\vskip .06in


\ref \key EW \by T.J. Enright, N.R. Wallach \pages 1-15
\paper Notes on homological algebra and representations of Lie algebras
\jour Duke Math. J.
{\bf 47} \yr 1980
\endref
\vskip .06in


\ref \key KZ \by A.W. Knapp, G. Zuckerman \pages 2178-2180
\paper Classification of irreducible tempered representations of semisimple Lie groups
\jour Proc. Nat. Acad. Sci. U.S.A. 
{\bf 73} \yr 1976
\endref
\vskip .06in


\ref \key Ko \by B. Kostant  \pages 329-387
\paper Lie algebra cohomology and the generalized Borel-Weil theorem
\jour Ann. of Math.
{\bf 74} \yr 1961
\endref
\vskip .06in


\ref \key L1\by R.P. Langlands \pages 253-257, A.M.S. Providence, R.I
\paper Dimension of spaces of automorphic forms, Algebraic Groups and Discontinuous Subgroups
\jour Proc. Sym. in Pure Math. 9
\yr 1966  
\endref
\vskip .06in


\ref \key L2\bysame \pages
\paper On the classification of irreducible representations of real algebraic groups
\newline
\paperinfo mimeographed notes
\jour Institute for Advanced Study \yr 1973
\endref
\vskip .06in


\ref \key M\by O. Mathieu\pages 537-592
\paper Classification of irreducible weight modules
\jour Ann. Inst. Fourier
{\bf 50} \yr 2000
\endref
\vskip .06in


\ref \key Mi\by I. Mirkovi\'c \pages 
\paper Classification of irreducible tempered representations of semisimple groups
\jour Ph.D. Dissertation, University of Utah, Salt Lake City \yr 1986
\endref
\vskip .06in


\ref \key PSZ \by I. Penkov, V. Serganova, G. Zuckerman \pages
\paper On the existence of $(\g,\k)$-modules of finite type
\paperinfo Duke Math. J., to appear
\endref
\vskip .06in


\ref \key PZ\by I. Penkov, G. Zuckerman \pages 311-326
\paper Generalized Harish-Chandra modules: a new direction in structure theory of representations
\jour Acta Applicandae Mathematicae 
{\bf 81} \yr 2004
\endref
\vskip .06in

\ref \key Sa \by G. Savin \pages 177-184
\paper Dual pair $PGL(2)\times G_2$ and $(\g_2, SL(3))$-modules
\jour IMRN  no. 4
\yr 1994
\endref
\vskip .06in

\ref \key Sc \by W. Schmid \pages
\paper Homogeneous complex manifolds and representations of semisimple Lie groups
\paperinfo Ph.D. Dissertation, University of California at Berkeley, 1967
\endref
\vskip .06in

\ref \key S\by J.P. Serre \pages
\paper Representations lin\'eares et espaces homogen\'es K\"ahlerians des groups de Lie compacts
\paperinfo Expos\'e 100, S\'eminarie Bourbaki, $6^e$ ann\'ee, (1953/54)
\endref
\vskip .06in

\ref \key V1\by D. Vogan \pages 1-60
\paper The algebraic structure of the representations of semisimple Lie groups I
\jour Ann. of Math. 
{\bf 109}  \yr 1979
\endref
\vskip .06in

\ref \key V2\bysame \pages
\paper Representations of Real Reductive Lie Groups
\jour Progress in Math., Birkhauser, Boston 15
\yr 1981
\endref
\vskip .06in

\ref \key  We \by C. Weibel\pages
\paper An Introdution to Homological Algebra
\paperinfo Cambridge Studies in Advanced Mathematics 38, Cambridge University Press, Cambridge
\yr 1994
\endref
\vskip .06in

\ref \key Z\by  G. Zuckerman\pages
\paper Construction of representations by derived functors
\paperinfo Lectures at the Institute for Advanced Study in Princeton, 1978 (unpublished notes)
\endref
